\documentclass[12pt,a4paper]{article}
\usepackage[cp1251]{inputenc}
\usepackage{amssymb}
\usepackage{amsmath}
\usepackage{amsthm}
\usepackage{euscript}

\newtheorem{theor}{Theorem}
\newtheorem{remark}{Remark}

\begin{document}

{\bf \centerline{On a sum of centered random variables with nonreducing variances. }}

\medskip
\centerline {I.~Shnurnikov}

\smallskip
\bigskip
{\footnotesize {\bf Abstract.}  Let $x=(x_1,\dots, x_n)$ and $\varepsilon=(\varepsilon_1,\dots, \varepsilon_n ) \in \mathbb{R}^n$. If $\|x\|_2=1$ and the coordinates $\varepsilon_i$ of $\varepsilon$  are independently distributed random variables with $Pr(\varepsilon_i=1) \ = \ Pr(\varepsilon_i=-1)= \frac 12$, then
$$
Pr(|\varepsilon^T x|\leq 1)\geq 0.36
$$
We give a new proof of this inequality which is weaker then the best known one $Pr(|\varepsilon^T x|\leq 1)\geq \frac 38$, proved by R.~Holzman and D.J.~Kleitman.
}

\bigskip
{\bf Introduction.}
Let us consider $x=(x_1,\dots, x_n)$ and $\varepsilon=(\varepsilon_1,\dots, \varepsilon_n ) \in \mathbb{R}^n$ such that $\|x\|_2=1$ and the coordinates $\varepsilon_i$ of $\varepsilon$  are independently distributed random variables with $Pr(\varepsilon_i=1) \ = \ Pr(\varepsilon_i=-1)= \frac 12$.

The problem is to give low bounds for $Pr(|\varepsilon^T x|\leq 1)$ for arbitrary $x$.
R.~Holzman and D.J.~Kleitman in \cite{Holzman}
proved that $Pr(|\varepsilon^T x|\leq 1)\geq \frac 38$ and conjectured that $Pr(|\varepsilon^T x|\leq 1)\geq \frac 12.$ For comparison the strong inequality $Pr(|\varepsilon^T x|< 1) \geq \frac 38$ is sharp. There are various bounds in \cite{Veraar}, see it also for references.

Developing robust optimization technology A.~Ben-Tal, A.~Nemirovski and C.~Roos have had to prove that $Pr(|\varepsilon^T x|\leq 1)$ is restricted from zero to proceed an effective bound of the "level of conservativeness"$,$ see \cite{Ben-Tal Nemirovsky} for details. It was proved in \cite{Ben-Tal Nemirovsky} that $Pr(|\varepsilon^T x|\leq 1)\geq \frac 13$ (result of R.~Holzman and D.J.~Kleitman was unknown yet).

In the present communication we shall show how to upgrade arguments of \cite{Ben-Tal Nemirovsky} to prove that $Pr(|\varepsilon^T x|\leq 1)\geq 0.36$.  We hope that this bound may be improved in a similar way up to more then $\frac 38$.

\begin{theor}
Let $x=(x_1,\dots, x_n)$ and $\varepsilon=(\varepsilon_1,\dots, \varepsilon_n ) \in \mathbb{R}^n$. If $\|x\|_2=1$ and the coordinates $\varepsilon_i$ of $\varepsilon$  are independently distributed random variables with $Pr(\varepsilon_i=1) \ = \ Pr(\varepsilon_i=-1)= \frac 12$, then
\begin{equation*}
Pr(|\varepsilon^T x|\leq 1)\geq 0.36
\end{equation*}
\end{theor}

\bigskip
\noindent
{\it Proof.}
Without loss of generality we may assume that
$$
x_1\geq x_2\geq \dots \geq x_n\geq 0.
$$
Let us consider two cases.

{\it Case 1}, $x_1+x_2>1.$ Let us denote
\begin{equation*}
r=r(\varepsilon)=\varepsilon_3x_3+\dots+\varepsilon_nx_n
\end{equation*}
The random quantity $r$ has symmetric distribution and hence
\begin{gather}
Pr(|\varepsilon^T x|\leq 1)=\frac 12Pr(-1\leq r+x_1+x_2\leq 1)+ \frac 12Pr(-1\leq r+x_1-x_2\leq 1)= \notag\\
=\frac 12 \left(Pr(-1-x_1-x_2 \leq r \leq 1-x_1-x_2)+Pr(-1-x_1+x_2 \leq r \leq 1-x_1+x_2)\right)=\notag\\
=\frac 12 \left( Pr(0\leq r \leq 1+x_1+x_2)+Pr(0<r \leq 1+x_1-x_2)\right)+\notag\\
+\frac 12 Pr(x_1+x_2-1 \leq r \leq 1-x_1+x_2)  \geq\notag\\
\geq \frac 14 \left( Pr(|r| \leq 1+x_1+x_2)+Pr(|r| \leq 1+x_1-x_2) \right). \label{Pr v summe}
\end{gather}

By the Tschebyshev inequality
\begin{equation}
\label{Tschebyshev r 2}
Pr(|r| \leq 1+x_1-x_2)>1-\frac{M r^2}{(1+x_1-x_2)^2} \geq \frac 12,
\end{equation}
since $x_1 \geq x_2$ and
\begin{equation*}
M r^2=1-x_1^2-x_2^2 \geq \frac 12
\end{equation*}
 for $x_1+x_2\geq 1.$

By the 4-th degree Tschebyshev inequality we have
\begin{equation*}
Pr(|r| \leq 1+x_1+x_2)>1-\frac{M r^4}{(1+x_1+x_2)^4}=1- \frac 1{(1+x_1+x_2)^4} \left( \sum_{i=3}^n x_i^4+6\sum_{3\leq i<j\leq n}x_i^2x_j^2\right).
\end{equation*}
Notice that $1+x_1+x_2 \geq 2$ and
\begin{equation*}
\sum_{i=3}^n x_i^4+6\sum_{3\leq i<j\leq n}x_i^2x_j^2<3(x_3^2+\dots+x_n^2)^2< \frac 34.
\end{equation*}
Hence
\begin{equation}
\label{Pr r Ch 4}
Pr(|r|\leq 1+x_1+x_2)>1-\frac 3{64}.
\end{equation}
So in the case $x_1+x_2>1$ inequalities (\ref{Pr v summe}), (\ref{Tschebyshev r 2}) and (\ref{Pr r Ch 4}) yield
\begin{equation*}
Pr(|\varepsilon^Tx|\leq 1)\geq \frac{93}{256} > 0.36
\end{equation*}

{\textit Case 2.} $x_1+x_2 \leq 1.$ Let us denote
\begin{equation*}
s_k=s_k(\varepsilon)=\varepsilon_1x_1+\dots+\varepsilon_k x_k
\end{equation*}
for $1\leq k \leq n$.
Let us define the following events
\begin{equation*}
A_k=\{\varepsilon\ : \ |s_j(\varepsilon)|\leq 1-x_{j+1}, \ j=1, \dots, k-1, \quad \text{and } \ |s_k(\varepsilon)|> 1-x_{k+1}\}
\end{equation*}
for $2 \leq k \leq n-1$ and
\begin{equation*}
A_n=\{\varepsilon\ : \ |s_j(\varepsilon)|\leq 1-x_{j+1},\ j=1, \dots n-1\}.
\end{equation*}
In the current case the events $A_2, \dots , A_n$ form a partition of probability space. In order to prove
\begin{equation*}
Pr(|s_n|\leq 1) \geq 0.36
\end{equation*}
 we shall get
\begin{equation*}
Pr(|s_n|\leq 1\ | \ A_k) \geq 0.36
\end{equation*}
provided $Pr(A_k)>0.$
Surely $Pr(|s_n|\leq 1\ | \ A_n)=1.$ If $\varepsilon \in A_k$ for $2\leq k \leq n-1,$ then $|s_k(\varepsilon)|\leq 1$ and
\begin{multline}
Pr(|s_n|\leq 1 \ | \ A_k)\geq Pr(0 \leq s_n-s_k \leq 2-x_{k+1})\geq \\  \frac 12 Pr(|s_n-s_k|\leq 2-x_{k+1})
 \geq \frac 12 \left(1  - \frac {M (s_n-s_k)^2}{(2-x_{k+1})^2} \right)
\label{Pr via M}
\end{multline}
by Tschebyshev inequality.
Note that
\begin{equation}
\label{M g ine}
M (s_n-s_k)^2=1-x_1^2-\dots-x_k^2\leq 1-kx_{k+1}^2.
\end{equation}
Let us define a function
\begin{equation*}
g_k(x)=\frac 12 \left(1 - \frac{1-kx^2}{(2-x)^2}\right).
\end{equation*}
It follows from (\ref{Pr via M}) and (\ref{M g ine}) that
\begin{equation}
\label{Pr > g_k}
Pr(|s_n|\leq 1 \ | \ A_k) \geq g_k(x_{k+1}).
\end{equation}

Since $\varepsilon \in A_k$, so
\begin{equation*}
1-x_{k+1}\leq x_1+x_2+\dots+x_k.
\end{equation*}
By the Cauchy inequality
\begin{equation*}
\frac{ x_1+x_2+\dots+x_k}k \leq \sqrt{\frac{x_1^2+\dots+x_k^2}{k}}=\sqrt{\frac{1-M(s_n-s_k)^2}{k}}.
\end{equation*}
Hence
\begin{equation}
\label{M h ine}
M(s_n-s_k)^2\leq 1- \frac {(1-x_{k+1})^2}{k}.
\end{equation}
Let us define a function
\begin{equation*}
h_k(x)=\frac 12 \left(1-\frac{1-\frac{(1-x)^2}{k}}{(2-x)^2}\right)
\end{equation*}
It follows from (\ref{Pr via M}) and (\ref{M h ine}) that
\begin{equation}
\label{Pr > h_k}
Pr(|s_n|\leq 1 \ | \ A_k) \geq h_k(x_{k+1}).
\end{equation}

Now we are going to analyse functions $g_k(x)$ and $h_k(x)$ for $x\in [0,1]$ and $k\geq 2$. It is easy to check that $g_k(x)$ increases for $x\in [\frac 1{2k},1]$ and that $h_k(x)$ decreases for $x\in [0,1]$. Moreover, $g_k(x)=h_k(x) \Leftrightarrow x=\frac 1{k+1}$ and
\begin{equation*}
\min_{x \in [0,1]} \max\left\{g_k(x),h_k(x)\right\}=g_k\left(\frac 1{k+1}\right).
\end{equation*}
This allows us to deduce from (\ref{Pr > g_k}) and (\ref{Pr > h_k}) the following
 \begin{equation*}
Pr(|s_n|\leq 1 \ | \ A_k) \geq \max\left\{g_k(x_{k+1}), \ h_k(x_{k+1})\right\}\geq g_k\left(\frac 1{k+1}\right)\geq g_2\left(\frac 13\right)=0.36.
\end{equation*}
Since the events $A_2, \dots , A_n$ form a partition of probability space, we finally have
 \begin{equation*}
Pr(|s_n|\leq 1) \geq 0.36.
\end{equation*}

\bigskip
\begin{remark}
One may use 4---th degree Tschebyshev inequality in the case 2 part of the proof instead of {\rm (\ref{Pr via M})}. It may lead to a better bound and I suppose, that $p \geq 0.4$ is possible to prove in a such way.
\end{remark}

\bigskip

{\bf Acknowledges}. Author is thankful to G.V.~Dilman, V.U.~Protasov and M. Veraar for helpful discussions and motivation.

\end{document}